\documentclass[12pt,leqno]{article}
\usepackage{amsfonts,amsthm,amsmath}
\theoremstyle{plain}

\newtheorem{thm}{Theorem}[section]

\newtheorem{lem}{Lemma}[section]

\theoremstyle{remark}
\newtheorem{rem}{Remark}[section]

\usepackage{setspace}
\newcommand{\ZZ}{\mathbb{Z}}
\newcommand{\RR}{\mathbb{R}}
\newcommand{\CC}{\mathbb{C}}

\DeclareMathOperator{\swe}{swe}

\begin{document}

\title{Nonexistence for extremal Type~II $\ZZ_{2k}$-Codes
}

\author{
Tsuyoshi Miezaki\thanks{
Division of Mathematics, 
Graduate School of Information Sciences, 
Tohoku University, 
6-3-09 Aramaki-Aza-Aoba, Aoba-ku, 
Sendai 980-8579, Japan. 
email: miezaki@math.is.tohoku.ac.jp
}
}
\date{}
\maketitle

\begin{abstract}
In this paper, we show that an extremal Type II $\ZZ_{2k}$-code
of length $n$ dose not exist for all sufficiently large $n$ when $k=2,3,4$.
\end{abstract}

%%%%%%%%%%%%%%%%%%%%%%%%%%%%%%%%%
\section{Introduction}
Let $\ZZ_{2k}\ (=\{0,1,2,\ldots,2k-1\})$ be the ring 
of integers modulo $2k$, where $k$ 
is a positive integer. 
We sometimes regard the elements of $\ZZ_{2k}$ as those of $\ZZ$. 
A $\ZZ_{2k}$-code $C$ of length $n$
(or a code $C$ of length $n$ over $\ZZ_{2k}$)
is a $\ZZ_{2k}$-submodule of $\ZZ_{2k}^n$.
A code $C$ is {\em self-dual} if $C=C^\perp$ where
the dual code $C^\perp$ of $C$ is defined as 
$C^\perp = \{ x \in \ZZ_{2k}^n \mid x \cdot y = 0$ for all $y \in C\}$
under the standard inner product $x \cdot y$. 
%% Two $\ZZ_{2k}$-codes $C$ and $C'$ are {\em equivalent} 
%% if there exists a monomial $(\pm 1, 0)$-matrix $P$ with 
%% $C' = C \cdot P = \{ x P\:|\: x \in C\}$.  
%% The automorphism group $\Aut(C)$ of $C$ is the group of all
%% monomial $(\pm 1, 0)$-matrices $P$ with
%% $C = C \cdot P$.
The Euclidean weight of a codeword 
$x=(x_1,x_2,\ldots,x_n)$ is 
$\sum_{i=1}^n \min\{x_i^2,(2k-x_i)^2\}$.
The minimum Euclidean weight $d_E(C)$ of $C$ is the smallest Euclidean
weight among all nonzero codewords of $C$.

A binary doubly even self-dual code is often called Type~II.
For $\ZZ_4$-codes, Type~II codes were first defined 
in~\cite{Z4-BSBM} as self-dual codes containing a 
$(\pm 1)$-vector and with the property that all
Euclidean weights are divisible by eight.
Then it was shown in~\cite{Z4-HSG} that,
more generally, the condition of containing a 
$(\pm 1)$-vector is redundant.
Type~II $\ZZ_{2k}$-codes was defined in 
\cite{BDHO} as a self-dual code with the property that all
Euclidean weights are divisible by $4k$.
It is known that a Type~II $\ZZ_{2k}$-code of length $n$ exists 
if and only if $n$ is divisible by eight.

In \cite{H-M}, we show the following theorem:
\begin{thm}[cf. \cite{H-M}]\label{main1}
Let $C$ be a Type~II $\ZZ_{2k}$-code of length $n$.
If $k \le 6$
then the minimum Euclidean weight $d_E(C)$ of $C$ 
is bounded by
\begin{equation}\label{Eq:B}
d_E(C) \le 4k \Big\lfloor \frac{n}{24} \Big\rfloor +4k.
\end{equation}
\end{thm}

\begin{rem}
The upper bound (\ref{Eq:B})
is  known for the cases $k=1$~\cite{MS73}
and $k=2$~\cite{Z4-BSBM}.
For $k \ge 3$, the bound (\ref{Eq:B}) is known under the
assumption that $\lfloor n/24 \rfloor\leq k-2$~\cite{BDHO}.
\end{rem}
In \cite{H-M}, we define that a Type~II $\ZZ_{2k}$-code meeting the bound (\ref{Eq:B})
with equality is {\em extremal} for $k \le 6$. 

The aim of this paper is to show the following theorem. 
\begin{thm}\label{main2}
For $k \le 4$, an extremal Type II $\ZZ_{2k}$-code
of length $n$ dose not exist for all sufficiently large $n$. 
\end{thm}

\begin{rem}
For the case $k=1$, the above result 
in Theorem \ref{main2} was shown in \cite{MS73}.
\end{rem}

%%%%%%%%%%%%%%%%%%%%
\section{Preliminaries}

An $n$-dimensional (Euclidean) lattice $\Lambda$ is a subset of
$\RR^{n}$ with the property that there exists a basis 
$\{e_{1}, e_2,\ldots, e_{n}\}$ of $\RR^{n}$ such that 
$\Lambda =\ZZ e_{1}\oplus \ZZ e_{2}\oplus \cdots \oplus\ZZ e_{n}$, i.e., 
$\Lambda$ consists of all integral linear combinations of the 
vectors $e_{1}, e_{2}, \ldots, e_{n}$. 
The dual lattice $\Lambda^*$ of $\Lambda$ is the lattice
$
\{x\in \RR^{n} \mid \langle x,y \rangle \in\ZZ \text{ for all }
y\in \Lambda\}
$, 
where $\langle x, y \rangle$ is the standard inner product.
%% A lattice $\Lambda$ is {\em integral} if
%% $\Lambda \subseteq \Lambda^{*}$.
A lattice with $\Lambda=\Lambda^{*}$
is called {\em unimodular}.
The norm of $x$ is $\langle x, x \rangle$.
A unimodular lattice with even norms is said to be {\em even}, 
otherwise {\em odd}.
%A unimodular lattice containing a vector of odd norm is  said to be
%{\em odd}.
An $n$-dimensional even unimodular lattice exists if and only
if $n \equiv 0 \pmod 8$, while an odd unimodular lattice
exists for every dimension.
The minimum norm $\min(\Lambda)$ of $\Lambda$ is the smallest
norm among all nonzero vectors of $\Lambda$.
For $\Lambda$ and a positive integer $m$, 
the shell $\Lambda_m$ of norm $m$ is defined as
$\{x\in \Lambda \mid \langle x,x \rangle=m \}$.

The theta series of $\Lambda$ 
is 
\begin{align*}
\Theta_{\Lambda}(z)=
\Theta_{\Lambda}(q)=
\sum_{x\in \Lambda}q^{\langle x, x \rangle}
=\sum_{m=0}^{\scriptstyle \infty}\vert \Lambda_{m}\vert q^{m},
\quad q=e^{\pi i z},\ \text{Im}(z)>0.
\end{align*}
%% The theta series $\Theta_{\Lambda}(q)$ of $\Lambda$ 
%% is the formal power series 
%% \begin{align*}
%% \Theta_{\Lambda}(q)=
%% \sum_{x\in \Lambda}q^{\langle x, x \rangle}
%% =\sum_{m=0}^{\scriptstyle \infty}\vert \Lambda_{m}\vert q^{m}. 
%% \end{align*}
For example, let $\Lambda$ be the $E_{8}$-lattice. Then, 
\begin{align*}
\Theta_{\Lambda}(q)=E_{4}(q)&=1+240\sum_{m=1}^{\infty}\sigma _{3}(m) q^{2m} \\
&=1 + 240 q^2 + 2160 q^4 + 6720 q^6 + 17520 q^8 +\cdots, 
\end{align*}
where $\sigma _{3}(m)$ is a divisor function $\sigma _{3}(m)=\sum_{0<d|m}d^3$. 

It is well-known that if $\Lambda$ is an $n$-dimensional 
even unimodular lattice, 
then $\Theta_{\Lambda}$ is a modular form of weight $n/2$ 
for the full modular group $SL_{2}(\ZZ)$ (see \cite{SPLAG}). 
For example, $E_{4}(q)$ is a modular form of weight $4$ 
for $SL_{2}(\ZZ)$. 
Moreover the following theorem is known (see \cite[Chap.~7]{SPLAG}).

\begin{thm}
\label{thm:Hecke}
If $\Lambda$ is an even unimodular lattice, then 
\[
\Theta _{\Lambda}(q)\in \CC[E_{4}(q), \Delta _{24}(q)],
\]
where 
$\Delta _{24}(q)=q^{2}\prod _{m=1}^{\scriptstyle \infty}
(1-q^{2m})^{24}$ 
which is a modular form of weight $12$ for $SL_{2}(\ZZ)$. 
\end{thm}

We now give a method to construct 
even unimodular lattices from Type II codes, which
is called Construction A \cite{BDHO}.
Let $\rho$ be a map from $\ZZ_{2k}$ to $\ZZ$ sending $0, 1, \ldots , k$ 
to $0, 1, \ldots , k$ and $k+1, \ldots , 2k-1$ to $1-k, \ldots , -1$, 
respectively. 
If $C$ is a  self-dual $\ZZ_{2k}$-code of length $n$, then 
the lattice 
\[
A_{2k}(C)=\frac{1}{\sqrt{2k}}\{\rho (C) +2k \ZZ^{n}\} 
\]
is an $n$-dimensional unimodular lattice, where 
\[
\rho (C)=\{(\rho (c_{1}), \ldots , \rho (c_{n}))\ 
\vert\ (c_{1}, \ldots , c_{n}) \in C\}. 
\]
The minimum norm of $A_{2k}(C)$ is $\min\{2k, d_{E}(C)/2k\}$.
Moreover, 
if $C$ is Type II, then the lattice $A_{2k}(C)$ 
is an even unimodular lattice.

The symmetrized weight enumerator of a $\ZZ_{2k}$-code
$C$ is 
\begin{align*}
\swe_{C}(x_{0}, x_{1}, \ldots , x_{k})=
\sum_{c\in C}x_{0}^{n_{0}(c)}x_{1}^{n_{1}(c)}\cdots 
x_{k-1}^{n_{k-1}(c)}x_{k}^{n_{k}(c)}, 
\end{align*}
where $n_{0}(c)$, $n_{1}(c)$, \ldots , 
$n_{k-1}(c)$, $n_{k}(c)$ are the number of 
$0, \pm 1, \ldots , \pm k-1, k$ components of $c$, respectively. 
Then the theta series of $A_{2k}(C)$ can be found by replacing 
$x_{1}$, $x_{2}$, $\ldots$, $x_{k}$ by 
\[
f_{0}=\sum_{x\in 2k\ZZ}q^{x^{2}/2k},\ 
f_{1}=\sum_{x\in 2k\ZZ +1}q^{x^{2}/2k},\ \ldots ,\ 
f_{k}=\sum_{x\in 2k\ZZ +k}q^{x^{2}/2k}. 
\]
respectively.
%
%%%%%%%%%%%%%%%%%%%%%%%%%%%%%
%\section{Review of Theorem~\ref{main1}}
%In this section, we review a proof of Theorem~\ref{main1}. 
%For the details, see \cite{H-M}. 
%%%%%%%%% proof of Theorem A %%%%%%%%%%
%{\it Proof of Theorem \ref{main1}.}
%Let $C$ be a Type~II $\ZZ_{2k}$-code of length $n$. Then 
Let $C$ be a Type II $\ZZ_{2k}$-code of length $n$. 
Then, the even unimodular lattice $A_{2k}(C)$ contains the sublattice 
$\Lambda _{0}=\sqrt{2k}\ZZ^{n}$ which has minimum norm $2k$.
We set $\Theta _{\Lambda_{0}}(q)=\theta _{0}$,
$n=8j$ and $j=3\mu +\nu$ ($\nu=0,1,2$), that is, $\mu=\lfloor n/24 \rfloor$.
%In this proof, 
We denote $E_4(q)$ and $\Delta_{24}(q)$ by
$E_{4}$ and $\Delta$, respectively.
By Theorem~\ref{thm:Hecke}, 
the theta series of $A_{2k}(C)$ can be written as 
\begin{align*}
\Theta _{A_{2k}(C)}(q)=\sum_{s=0}^{\mu}a_{s}E_{4}^{j-3s}\Delta ^{s}
=\sum_{r\geq 0}\vert {A_{2k}(C)}_{r}\vert q^{r}
=\theta_{0}+\sum_{r\geq 1}\beta _{r}q^{r}. 
\end{align*}
%Suppose that $d_{E}(C)\geq 4k\lfloor n/24 \rfloor +4k$. 
Let $C$ be an extremal Type II $\ZZ_{2k}$-code for 
$1\leq k\leq 6$, namely, 
$d_{E}(C)= 4k(\mu +1)$. 
We remark that a codeword of Euclidean weight $4km$
gives a vector of norm $2m$ in $A_{2k}(C)$. 
Then we choose the $a_{0}, a_1, \ldots , a_{\mu}$ so that 
\begin{align*}
\Theta _{A_{2k}(C)}(q)=\theta_{0}+\sum_{r\geq 2(\mu+1)}\beta ^{\ast}_{r}q^{r}. 
\end{align*}

Here, we set $b_{2s}$ as 
$E_{4}^{-j}\theta_{0}=\sum_{s=0}^{\infty}b_{2s}(\Delta/E_{4}^{3})^s$. 
That is, $\theta_{0}=\sum_{s=0}^{\infty}b_{2s}E_{4}^{j-3s}\Delta^s$. 
Then 
\begin{align*}
\sum_{s=0}^{\mu}a_{s}E_{4}^{j-3s}\Delta ^{s}
=\Theta _{A_{2k}(C)}(q)
=\sum_{s=0}^{\infty}b_{2s}E_{4}^{j-3s}\Delta^s
+\sum_{r\geq 2(\mu+1)}\beta ^{\ast}_{r}q^{r}. 
\end{align*}
Comparing the coefficients of $q^{i}$ $(0\leq i\leq 2\mu)$, we get 
$a_{s}=b_{2s}$ $(0\leq s\leq \mu)$. 
Hence we have 
\begin{align}
-\sum_{r\geq (\mu+1)}b_{2r}E_{4}^{j-3r}\Delta^r
=\sum_{r\geq 2(\mu+1)}\beta^{\ast}_{r}q^{r}. \label{eqn:b}
\end{align}
In (\ref{eqn:b}), comparing the coefficients of $q^{2(\mu+1)}$ 
and $q^{2(\mu+2)}$, 
we have
\begin{align}
\left\{
\begin{array}{ll}
\beta^{\ast}_{2(\mu +1)}&=-b_{2(\mu +1)}, \\
\beta^{\ast}_{2(\mu +2)}&=-b_{2(\mu +2)}+b_{2(\mu +1)}(24\mu -240\nu +744). 
\end{array}
\right. \label{eqns:coef}
\end{align}
All the series are in $q^{2}=t$, 
and B\"urman's formula {\cite[page 128]{WW}} shows that 
\begin{align}\label{eqn:B}
b_{2s}=\frac{1}{s!}\frac{d^{s-1}}{dt^{s-1}}
\left(
\left(
\frac{d}{dt}(E_{4}^{-j}\theta _{0})\right)
(t E_{4}^{3}/\Delta )^{s}\right)_{\{t=0\}}. 
\end{align}
In \cite{H-M}, we show that 
\begin{align}\label{inequ:main1}
\beta^{\ast}_{2(\mu +1)}>0 
\end{align}
and we remark that the inequality (\ref{inequ:main1}) 
is a crucial part of the proof of Theorem \ref{main1}. 

Finally, we quote the two theorems needed later: 
\begin{thm}[cf.~{\cite[page 18, Theorem 1.64]{Web}}]\label{thm:eta1}
Let $\eta(z)=t^{1/24}\prod_{m=1}^{\infty}(1-t^{m})$ 
be the Dedekind eta function, where $t=e^{2\pi iz}$, 
the same for several places and 
$Im(z)>0$. 
If $f(z)=\prod_{\delta \vert N}\eta(\delta z)^{r_{\delta}}$ 
with $k=(1/2)\sum_{\delta\vert N}r_{\delta}\in\ZZ$, 
with the additional properties that 
\[
\sum_{\delta\vert N}\delta r_{\delta}\equiv 0\pmod{24}
\]
and 
\[
\sum_{\delta\vert N}\frac{N}{\delta}r_{\delta}\equiv 0\pmod{24}, 
\]
then $f(z)$ satisfies 
\[
f\left(\frac{az+b}{cz+d}\right)=\chi(d)(cz+d)^{k}f(z)
\]
for every 
${\small\begin{pmatrix}
a&b\\
c&d
\end{pmatrix}}\in\Gamma_{0}(N)$. 
Here the character $\chi$ is defined by $\chi(d):=\left(\frac{(-1)^{k}s}{d}\right)$, 
where $\Big(\frac{\cdot}{\cdot}\Big)$ is the usual 
Jacobi symbol and $s:=\prod_{\delta \vert N}\delta ^{r_{\delta}}$. 
\end{thm}
\begin{thm}[cf.~{\cite[page 18, Theorem 1.65]{Web}}]\label{thm:eta2}
Let $c$, $d$ and $N$ be positive integers with 
$d\vert N$ and $\gcd(c,d)=1$. 
If $f(z)=\prod_{\delta \vert N}\eta(\delta z)^{r_{\delta}}$ satisfying 
the conditions of Theorem \ref{thm:eta1} for $N$, 
then the order of vanishing of $f(z)$ at the cusp $c/d$ is 
\[
\frac{N}{24}\sum_{\delta\vert N}\frac{\gcd(d,\delta)^{2}r_{\delta}}
{\gcd(d,\frac{N}{d})d\delta}. 
\]
\end{thm}
%%%%%%%%%%%%%%%%%%
\section{Proof of Theorem \ref{main2}}

%%%%%%%%% proof of Theorem B %%%%%%%%%%
In this section, we give a proof of Theorem~\ref{main2}.
Our proof is an analogue of that 
of~\cite{MOS75}. 
Before we give the proof of Theorem \ref{main2},
we give two lemmas. 
First, we quote the following lemma from \cite{MOS75}. 
In \cite{I}, Ibukiyama remarked that 
in ~\cite[Lemma 1]{MOS75} 
$2\pi$ (p.~70, l.~$-1$)
should be $(2\pi)^{1/2}$.

\begin{lem}[{\cite[Lemma 1]{MOS75}},~{\cite[Theorem 12]{I}}]\label{lem:MOS}
Suppose that $G(q)$, $H(q)$ are analytic inside the circle 
$\vert q\vert =1$ and satisfy{\rm :} 
\begin{align*}
&({\rm i})\ H(q)=\sum_{s=0}^{\infty} H_{s}q^{s}\ with\ H_{0}>0,\ 
H_{1}>0,\ and\ H_{s}\geq 0\ for\ s\geq 2, \\
&({\rm ii})\ if\ F(y)=e^{2\pi y}H(e^{-2\pi y}),
\ then\ F^{\prime }(y) =0\ has\ a\ solution\ 
y=y_{0} \ in\ the\ range\\ 
&\hspace{17pt}y>0,\ with\ F(y_{0})=c_{1}>0, 
F^{\prime \prime }(y_{0})/F(y_{0})=c_{2}>0, G(e^{-2\pi y_{0}})\neq 0. 
\end{align*}
Then $\beta_{r}$, the coefficient of $q^{r}$ in $G(q)H(q)^{r}$, 
satisfies 
\begin{align*}
\beta_{r}\sim \frac{(2\pi)^{1/2} }{(rc_{2})^{1/2}}
G(e^{-2\pi y_{0}})c_{1}^{r},\ as\ r\rightarrow \infty . 
\end{align*}
\end{lem}

Second, we show the following lemma: 
\begin{lem}\label{lem:nonzero}
We set $t=q^2=e^{2\pi iz}$ and $f_{0}(k,t)=\sum_{x\in \ZZ}t^{kx^{2}}$. 
Let $Z(k,t):=[f_0(k, t)^8, E_4(t)]/4=f_0(k, t)^8E_4(t)^{\prime}
-(f_0(k, t)^8)^{\prime}E_4(t)$, 
where $[\ ,\ ]$ is the Rankin-Cohen bracket and $f(t)^{\prime} 
= t (df / dt)$. Then, for $1\leq k \leq 4$ and a positive real number $y$, 
$Z(k,e^{(-2\pi y)})\neq 0$. 
\end{lem}
\begin{proof}
Let $f$ (resp. $g$) be a modular form of weight $k$ (resp. $\ell$) 
for a group $\Gamma$. Then, $[f,g]:=kfg^{\prime}-\ell f^{\prime}g$ is 
a modular form of weight $k+\ell+2$ for $\Gamma$ \cite[page 53]{1-2-3}. 

We remark that $f_{0}(1,t)$ is a modular form of weight $1/2$ 
for $\Gamma_{0}(4)$ \cite[page 12]{Web}. Therefore, 
$f_{0}(1,t)^4$ is a modular form of weight $2$ for $\Gamma_{0}(4)$. 
Moreover, $f_{0}(k,t)^4$ is a modular form of weight $2$ 
for $\Gamma_{0}(4k)$ \cite[page 28, Proposition 2.22]{Web}. 

\begin{itemize}
\item 
The case of $k=1$: 

We remark that $Z(1,t)\in \Gamma_0(4)$ and 
define the functions: 
\begin{align*}
\left\{
\begin{array}{l}
\Delta_4^{\infty}(t)=\eta^8(4z)/\eta^4(2z),\\
\Delta_4^{0}(t)=\eta^8(z)/\eta^4(2z),\\
J_4(t)=\Delta_4^{0}(t)/\Delta_4^{\infty}(t), 
\end{array}
\right.
\end{align*}
Note that $J_4(t)$ is an isomorphism from a 
fundamental domain of $\Gamma_0(4)$ to 
the Riemann sphere $\CC \cup \{\infty\}$ 
and a generator of 
the function field of $\mathbb{H}^{*} /\Gamma_0(4)$, 
where $\mathbb{H}$ be the upper half plane and 
$\mathbb{H}^{*} /\Gamma_0(4)$ is a compactification of 
$\mathbb{H}/\Gamma_0(4)$ \cite[page 407]{B}, 
\cite[page 16]{BKM}. 
Then, we have the following equality: 
\begin{align}\label{equ:nonzero4}
\frac{Z(1,t)}{\Delta_4^{\infty}(t)^5}
=224 X^4 + 11264 X^3 + 188416 X^2 + 1048576 X, 
\end{align}
where $X:=J_4(t)$. 
It is easy to check that 
there are no positive real roots of the right-hand side (\ref{equ:nonzero4}). 
Here, we remark that $J_4(e^{(2\pi iz)})$ takes a real 
on the imaginary axis. 
Using Theorem \ref{thm:eta1} and \ref{thm:eta2}, 
we have $\Delta_4^{\infty}(e^{(2\pi i 0)})\neq 0$ and 
$\Delta_4^{0}(e^{(2\pi i 0)})=0$, namely 
$J_4(e^{(2\pi i 0)})=0$. 
Therefore, the values of the $J_4(t)$ on the imaginary axis 
are positive real numbers and we have $Z(1,t)\neq 0$ on the imaginary axis. 
\end{itemize}
The other cases can be proved similarly. 
We only mention the functions which could 
be used for the proofs of the cases $k=2$, $3$ and $4$. 
\begin{itemize}
\item 
The case of $k=2$: \\
\begin{align*}
\left\{
\begin{array}{l}
\Delta_8^{\infty}(t)=\eta^4(8z)/\eta^2(4z),\\
\Delta_8^{0}(t)=\eta^4(z)/\eta^2(2z),\\
J_8(t)=\Delta_8^{0}(t)/\Delta_8^{\infty}(t), 
\end{array}
\right.
\end{align*}
where $J_8(t)$ is Hauptmodul for type ``$8-$" \cite[page 331]{CN}. 
\begin{align*}
Z(2,t)/\Delta_8^{\infty}(t)^{10}=&240 X^9+ 12928 X^8+ 283136 X^7+ 3358720 X^6\\
 &+ 
 23883776 X^5 + 105086976 X^4+ 281018368 X^3\nonumber\\
&+ 419430400 X^2  + 268435456 X  
\nonumber 
\end{align*}
where $X:=J_8(t)$. 

\item 
The case of $k=3$: \\
\begin{align*}
\left\{
\begin{array}{l}
\Delta_{12}^{\infty}(t)=\eta(2z)\eta^{-2}(4z)\eta^{-3}(6z)\eta^6(12z),\\
\Delta_{12}^{0}(t)=\eta^6(z)\eta^{-3}(2z)\eta^{-2}(3z)\eta(6z),\\
J_{12}(t)=(\Delta_{12}^{0}(t)/\Delta_{12}^{\infty}(t))^{1/2},
\end{array}
\right.
\end{align*}
where $J_{12}(t)$ is Hauptmodul for type ``$12-$" \cite[page 331]{CN}. 
\begin{align*}
Z(3,t)=
&240 X^{19}+ 18000 X^{18}+ 616032 X^{17}+12860832 X^{16}\\
&+ 184227840 X^{15}+ 1927623168 X^{14}+15293558784 X^{13}\nonumber\\
&+ 94189206528 X^{12}+ 456914313216 X^{11} + 
 1760257683456 X^{10}\nonumber \\ 
&+ 5401844490240 X^9+ 13181394788352 X^8+ 
25400510447616 X^7 \nonumber \\  
&+ 38149727846400 X^6+ 43699899727872 X^5 + 
 36857648775168 X^4 \nonumber\\  
&+ 21565588635648 X^3+ 7815347306496 X^2+1320903770112 X \nonumber
\end{align*}
where $X:=J_{12}(t)$.

\item 
The case of $k=4$: \\
\begin{align*}
\left\{
\begin{array}{l}
\Delta_{16}^{\infty}(t)=\eta(16z)^2/\eta(8z),\\
\Delta_{16}^{0}(t)=\eta^2(z)/\eta(2z),\\
J_{16}(t)=\Delta_{16}^{0}(t)/\Delta_{16}^{\infty}(t),
\end{array}
\right.
\end{align*}
where $J_{16}(t)$ is Hauptmodul for type ``$16-$" \cite[page 331]{CN}. 
\begin{align*}
Z(3,t)=
&240 X^{19}+ 13440 X^{18}+ 339840 X^{17}+5259776 X^{16}\\
&+56422912 X^{15}+448143360 X^{14}+2741043200 X^{13}\nonumber\\
&+13230211072 X^{12}+51153629184 X^{11} +159735971840 X^{10}\nonumber \\ 
&+403939164160 X^9+825259589632 X^8+1351740293120 X^7 \nonumber \\  
&+1750333390848 X^6+1751407132672 X^5 +1305938493440 X^4\nonumber\\  
&+682899800064 X^3+223338299392 X^2+34359738368 X \nonumber
\end{align*}
where $X:=J_{16}(t)$. 
\end{itemize}
\end{proof}

%% We begin the proof of Theorem \ref{main2}.

\bigskip
\noindent
{\it Proof of Theorem \ref{main2}.}
Using the equation (\ref{eqn:B}) and 
the fact that $\theta _{0}=\theta_{1}^{j}$ where $\theta _{1}$ 
is the theta series of the lattice $(2k\ZZ)^{8}/\sqrt{2k}$, 
we have 
\begin{align*}
b_{2s}=\frac{-j}{s!}\frac{d^{s-1}}{dt^{s-1}}\left(
E_{4}^{3s-j-1}\theta _{1}^{j-1}(\theta_{1}E_{4}^{\prime}
-\theta_{1}^{\prime}E_{4})(t/\Delta)^{s}\right)_{\{t=0\}}, 
\end{align*}
where $f^{\prime}$ is the derivation of $f$ with respect to $t=q^{2}$. 

We show that
$\beta^{\ast}_{2(\mu +2)}<0$ for sufficiently large $n$. 
When we set $h(t)=\prod _{r=1}^{\infty}(1-t^{r})^{-24}$, 
we have 
\begin{align*}
b_{2(\mu +1)}&=\frac{-j}{(\mu +1)!}\frac{d^{\mu}}{dt^{\mu}}\left(
E_{4}^{2-\nu}\theta _{1}^{j-1}(\theta_{1}E_{4}^{\prime}
-\theta_{1}^{\prime}E_{4})(h(q))^{\mu +1}\right)_{\{t=0\}}, \\
b_{2(\mu +2)}&=\frac{-j}{(\mu +2)!}\frac{d^{\mu +1}}{dt^{\mu +1}}\left(
E_{4}^{5-\nu}\theta _{1}^{j-1}(\theta_{1}E_{4}^{\prime}
-\theta_{1}^{\prime}E_{4})(h(q))^{\mu +2}\right)_{\{t=0\}}. 
\end{align*}

We show that 
$\vert b_{2(\mu +2)}/b_{2(\mu +1)} \vert$ 
is bounded, which implies that $\beta^{\ast}_{2(\mu +2)}<0$ 
as $n\rightarrow \infty$ since the equations (\ref{eqns:coef}) and 
the inequality (\ref{inequ:main1}) hold. 

We now apply Lemma \ref{lem:MOS} with 
$G(t)=G_{1}(t)=E_{4}^{2-\nu}\theta _{1}^{j-1}(\theta_{1}E_{4}^{\prime}
-\theta_{1}^{\prime}E_{4})h(t)$ and $H(t)=h(t)$. 
Then, as is shown in \cite{MOS75}, and using Lemma \ref{lem:nonzero}, 
the hypotheses (i) and (ii) in Lemma \ref{lem:MOS} are satisfied. So, 
\begin{align*}
b_{2(\mu +1)}\sim -(2\pi)^{1/2} jc_{2}^{-1/2}\mu ^{-3/2}
G_{1}(e^{-2\pi y_{0}})c_{1}^{\mu},\ as\ r\rightarrow \infty . 
\end{align*}
where 
%$c_{1}(=69.1\ldots)$
$c_{1}$ and $c_{2}$ are constants. 
Similarly with 
$G(q)=G_{2}(q)=E_{4}^{5-\nu}\theta _{1}^{j-1}(\theta_{1}E_{4}^{\prime}
-\theta_{1}^{\prime}E_{4})h(q)$ and $H(q)=h(q)$. 
\begin{align*}
b_{2(\mu +2)}\sim -(2\pi)^{1/2} jc_{2}^{-1/2}\mu ^{-3/2}
G_{2}(e^{-2\pi y_{0}})c_{1}^{\mu +1},\ as\ r\rightarrow \infty . 
\end{align*}
Hence $\vert b_{2(\mu +2)}/b_{2(\mu +1)} \vert$ 
is bounded (In fact, it approaches a limit of about $1.64\times 10^{5}$ 
as $\mu \rightarrow \infty$). 
%\end{proof}
\qed \bigskip

\begin{rem}
Using the equations (\ref{eqns:coef}), the coefficient 
$\beta^{\ast}_{2(\mu +2)}$ first goes negative 
when $n$ is about $1.64\times 10^{5}$. 
%Namely, for $k\le 4$, an extremal 
%Type II $\ZZ_{2k}$-code of length $n$ does not exist 
%when $\beta^{\ast}_{2(\mu +2)}$ is negative. 
%for $n\geq 5200$. 
%Hence there are finitely many extremal Type II $\ZZ_{2k}$-codes
%for $k \le 4$.
\end{rem}

\begin{rem}
For $k=5$ and $6$, we could not show $G(e^{-2\pi y_{0}})\neq 0$ 
in the hypothesis (ii) in Lemma \ref{lem:MOS}. 
The method of Lemma \ref{lem:nonzero} does not 
work because there are no Hauptmoduls for the groups 
$\Gamma_{0}(20)$ and $\Gamma_{0}(24)$ 
since the groups are not genus zero. 
%If $G(e^{-2\pi y_{0}})\neq 0$ holds for $k=5$ and $6$, 
%then, we have the same conclusion of Theorem \ref{main2} that 
%an extremal Type II $\ZZ_{2k}$-code
%of length $n$ dose not exist for all sufficiently large $n$. 
\end{rem}

\bigskip
\noindent
{\bf Acknowledgment.}
The author would like to thank Masaaki Harada, 
Kenichiro Tanabe and Junichi Shigezumi for useful discussions. 
The author would also like to thank 
the referee for informing us the reference \cite{I} 
and its helpful comments. 
The author was supported by JSPS Research Fellowship.

%%%%%%%%%%%%%%%%%  References  %%%%%%%%%%%%%%%%%%%%%%%%

\end{document}